\documentclass[10pt]{amsart}
\usepackage{amsxtra,latexsym,amssymb}
\usepackage{epsfig,rotate,amsthm}

\def\qed{{\hfill $\Box$}}

\def\Z{{\mathbb Z}}
\def\C{{\mathbb C}}

\theoremstyle{theorem}
\newtheorem{thm}{Theorem}[section]
\newtheorem{cor}{Corollary}[section]
\newtheorem{prop}{Proposition}[section]

\newtheorem{lem}{Lemma}[section]
\theoremstyle{definition}
\newtheorem{defn}{Definition}[section]
\theoremstyle{remark}

\newtheorem{rem}{\bf Remark}[section]

\begin{document}
\title[Whittaker model of $R(f)$]{On Whittaker modules over a class of
  algebras similar to $U(sl_{2})$}
\author{Xin Tang}
\address{Department of Mathematics \& Computer Science\\
Fayetteville State University\\
Fayetteville, NC 28301}
\email{xtang@uncfsu.edu}
\keywords{Algebras similar to $U(sl_{2})$, Whittaker model, Whittaker modules}
\thanks{}
\date{Oct 2006}
\subjclass[2000]{Primary 17B55, 20G; Secondary 17B50}
\maketitle
\begin{abstract}
Motivated by the study of invariant rings of finite groups on the
first Weyl algebras $A_{1}$ (\cite{AHV}) and finding interesting 
families of new noetherian rings, a class of algebras similar to 
$U(sl_{2})$ were introduced and studied by Smith in \cite{S}. Since 
the introduction of these algebras, research efforts have been 
focused on understanding their weight modules, and many important 
results were already obtained in \cite{S} and \cite{Ku}. But it seems 
that not much has been done on the part of nonweight modules. In 
this note, we generalize Kostant's results in \cite{K} on the 
Whittaker model for the universal enveloping algebras $U(\frak g)$ 
of finite dimensional semisimple Lie algebras $\frak g$ to Smith's 
algebras. As a result, a complete classification of irreducible 
Whittaker modules (which are definitely infinite dimensional) for 
Smith's algebras is obtained, and the submodule structure of any 
Whittaker module is also explicitly described.
\end{abstract}
%\normalsize
\section*{Introduction} \label{Intr}
Motivated by the study of the invariant rings of finite groups on the
first Weyl algebra $A_{1}$ and other important things, an interesting 
class of algebras $R(f)$ similar to $U(sl_{2})$ were introduced and 
studied by Smith in \cite{S}. Each algebra $R(f)$ is generated by 
three generators $E,\,F,\,H$ subject to the relations $EF-FE=f(H), \, 
HE-EH=E,$ and $ HF-FH=-F$, where $f$ is a polynomial in $H$. 
These algebras serve as a subclass of Witten's 7-parameter
deformations of $U(sl_{2})$ as studied in \cite{Ku}. As their 
name indicates, these algebras share a lot of similar properties 
with $U(sl_{2})$. The ring theoretic properties and the highest 
weight modules were first investigated in detail in \cite{S}. These 
algebras are somewhat commutative noetherian domain, and have the 
GK-dimension $3$ (\cite{S}). The center $Z(R)$ of $R(f)$ is also
proved to be isomorphic to the polynomial ring in one variable. 
The primitive ideals are classified by Smith (\cite{S}). Furthermore, 
a similar theory of highest weight modules and the category $\frak O$ 
is also constructed for these algebras by Smith (\cite{S}). In
particular, for some special parameters $f$, all finite dimensional 
representations of $R(f)$ are semisimple. For more details, we refer 
the reader to \cite{S}. These algebras have also been further studied 
in \cite{HS} and \cite{Ku} from the points of views of both ring 
theoretic properties and representation theory.

Since the introduction of these algebras, a lot of research efforts 
have been focused on trying to understand their weight modules 
(\cite{S}, \cite{Ku}). But it seems to us that not much has been done
for the part of nonweight modules. So it might be useful to present 
some specific constructions for nonweight irreducible modules over 
these algebras. In this paper, we are able to work out such a
possibility by generalizing Kostant's results on the Whittaker model 
for the universal enveloping algebras $U(\frak g)$ of finite
dimensional semisimple Lie algebras $\frak g$ to Smith's algebras
$R(f)$. As an application, we obtain a complete classification of all 
irreducible Whittaker modules, and the submodule structure of any
Whittaker module is also completely determined.

The initial investigation of the Whittaker model and hence Whittaker
modules for semisimple Lie algebras was started by Kostant in the 
seminal paper \cite{K}. The study of Whittaker modules is closely 
related the Whittaker equations and has nice applications in the 
theory of Toda lattice. For a nonsingular character of the nilpotent 
subalgebra $\frak n^{+}$ of $\frak g$, Kostant introduced the
Whittaker model of the center $Z(\frak g)$ of $U(\frak g)$ for 
finite dimensional semisimple Lie algebras $\frak g$. Whittaker 
model was used to study the structure of Whittaker modules over 
$U(\frak g)$ and several important structure theorems were proved by
Kontant for Whittaker modules in \cite{K}. Note that Whittaker
modules are very similar to Verma modules. But Whittaker modules 
have a special feature in that they are irreducible if and only if
they admit a central character, which is is not always the case for 
Verma modules. The Whittaker model was later on generalized and
studied for singular characters of $\frak n^{+}$ by Lnych in his Ph.D. 
thesis \cite{L}. Other similar works on this subject also appeared 
in \cite{MS1} and \cite{MS2}. As a matter of fact, Verma 
modules and Whittaker modules are two extreme cases of generalized 
Whittaker modules ( \cite{L}, \cite{MS1} and \cite{MS2}). Furthermore,
generalized Whittaker modules are mapped to holonomic $D-$modules on the flag 
variety of $\frak g$ via the Beilinson-Bernstein localization
(\cite{BB}). Based on this observation, a geometric study of Whittaker 
modules for finite dimensional semisimple Lie algebras was carried out 
in \cite{MS1} and \cite{MS2}.

In addition, a quantum analogoue of the Whittaker model has been constructed 
by Sevoastyanov for the topological version $U_{h}(\frak g)$ of 
quantized enveloping algebras by using their realizations via 
Coxeter elements in \cite{S1}. The major difficutly of a direct 
generalization of Kostant's results to the quantized case lies 
in the fact that there is no nonsingular character for the positive 
part of the quantized enveloping algebras because of the quantized 
Serre relations. To resolve this issue, he has to turn to the 
topological version of 
quantized enveloping algebras which has different realizations 
admitting nonsingular characters for the positive part. In the case 
of $\frak g=sl_{2}$, the situation is slightly different, since 
the quantized Serre relations are vaccum. Thus a direct generalization 
of Kostant's approach should work. And this has recently been worked 
out by Ondrus in \cite{O}. We have to admit that it is just a pure
luck that a similar pattern works for Smith's algebras.

Now let us mention a bit about the organization of this paper. In
Section 1, we recall the definition of Smith's algebras and some basic
results on their properties. In Section 2, we construct the Whittaker
model of the center $Z$ of $R(f)$, and classify all irreducible Whattaker
modules. In Section 3, we investigate the submodule structure of any
Whittaker module. Throughout this paper, the base field will be
assumed to be $\mathbb{C}$, though the results hold over any
algebraically closed field of characteristic zero.

\section{Algebras similar to $U(sl_{2})$}\label{S:algebras}
In this section, we recall the definition and basic properties of the
algebras $R(f)$ similar to $U(sl_{2})$ as introduced by Smith in \cite{S}.
\begin{defn}
(See \cite{S}) Let $f$ be a polynomial in $H$, the algebra $R(f)$ is
defined to be the $\mathbb{C}-$algebra generated by $E,F,H$ subject to 
the following relations:
\[
EF-FE=f(H),\,HE-EH=E,\,HF-FH=-F.
\]
and $R(f)$ is called an algebra similar to $U(sl_{2})$. We will sometimes
denote it by $R$ in short.
\end{defn}
\begin{prop}
(See \cite{S}) $R(f)$ has $GK-$dimenion $3$.
\end{prop}
\qed 
\begin{prop}
(See \cite{S}) $R(f)$ is a somewhat commutative algebra.
\end{prop}
\qed

\begin{cor}
If $V$ is a simple module, then every element of $End_{R(f)}(V)$ is a scalar.
\end{cor}
{Proof:} This follows from Quillen's Lemma and the fact $R$
is a somewhat commutative algebra (\cite{S}). For more detail, we
refer the reader to \cite{S}.
\qed

Let $R(E)$ denote the subalgebra of $R(f)$ generated by $E$, and $R(F,H)$ the 
subalgebra generated by $F,H$. Then we have 
\begin{prop}
$R(F,H)$ is isomorphic to  the enveloping algebra of the two 
dimensional nonabelian Lie algebra.
\end{prop}
{\bf Proof:} The proof is obvious. \qed

From \cite{S}, we have the following fact:
\begin{prop}
For any polynomial $f(H)$, there exist another polynomial $u(H)$ such that $f=\frac{1}{2}\Delta(u)=\frac{1}{2}(u(H+1)-u(H))$.
\end{prop}

\qed

In addition, $R(f)$ has a Casimir element $\Omega$ which is defined as
$\Omega=2FE+u(H+1)$ 
and a simple caculation shows the following:
\begin{prop}
(See \cite{S}) The center $Z(R)$ of $R$ is a polynomial ring 
generated by one variable $\Omega=2FE+u(H+1)$ over $\mathbb{C}$, 
where $u$ is the polynomial such that $f(H)=1/2\Delta(u)$.
\end{prop}
\qed

\section{The Whittaker Model for the center $Z(R)$ of $R(f)$}

In this seciton, we work out the Whittaker model for the center
$Z(R(f))$ of $R(f)$, and use it to study Whittaker modules over
$R(f)$. We obtain similar results as in \cite{K}. In fact, We will 
closely follow the formulation in \cite{K} with some slight modifications.
\begin{defn}
An algebra homomorphism $\eta\colon R(E)\longrightarrow \mathbb{C}$ is 
called a nonsingular character of $R(E)$ if $\eta(E)\neq 0$.
\end{defn}
\begin{defn}
Let $V$ be an $R-$module, a vector $v\in V$ is called a Whittaker
vector of type $\eta$ if $E$ acts on it through a nonsingular
character $\eta$, i.e., $Ev=\eta(E)v$. If $V=Rv$, then we call $V$ 
a Whittaker module of type $\eta$, and $v$ is called a cyclic 
Whittaker vector of type $\eta$. 
\end{defn}
From now on, we fix such a nonsingular character $\eta$ of $R(E)$. 
The following proposition follows from the triangular decomposition 
of $R$ in \cite{S}:
\begin{prop}
$R$ is isomorphic to $R(F,H)\otimes R(E)$ as vector spaces and $R$ 
is a free module over $R(E)$.
\end{prop}
\qed

Let $\eta \colon R(E)\longrightarrow \C$ be the fixed nonsingular
character of $R(E)$, and we denote by $R_{\eta}(E)$ the kernel of the 
character $\eta$. Then we have
\begin{prop}
$R(E)=\C \oplus R_{\eta}(E)$. Thus $R\cong R(F,H)\oplus RR_{\eta}(E)$.
\end{prop}
{\bf Proof:} Since $R(E)=\C \oplus R_{\eta}(E)$ and $R=R(F,H)\otimes
( \C \oplus R_{\eta}(E))$, so $R \cong R(F,H) \oplus RR_{\eta}(E)$.
\qed

Now we can define a projecton $\pi\colon R\longrightarrow R(F,H)$ 
from $R$ into $R(F,H)$ by taking the $R(F,H)$ component of any $u\in
R$. We denote the image $\pi(u)$ of $u$ by $u^{\eta}$ for short.
\begin{lem}
If $v\in Z(R)$, then we have $u^{\eta}v^{\eta}=(uv)^{\eta}$.
\end{lem}
{\bf Proof:} Let $u,v \in Z(R)$, then we have
\[ uv-u^{\eta}v^{\eta}=(u-u^{\eta})v+u^{\eta}(v-v^{\eta})\\
=v(u-u^{\eta})+u^{\eta}(v-v^{\eta})\in RR_{\eta}(E).
\]
So $(uv)^{\eta}=u^{\eta}v^{\eta}$.

\qed

\begin{prop}
$\pi\colon Z(R)\longrightarrow R(F,H)$ is an algebra isomorphism of
$Z(R)$ onto its image $W(F,H)$ in $R(F,H)$.
\end{prop}
{\bf Proof:} It follows from that above lemma that $\pi$ 
is a homomorphism of algebras. Since $Z(R)=\C[\Omega]$ 
and $\pi(\Omega)=2\eta(E)F+u(H+1)$ which is not zero in $W(F,H)$,
so $\pi$ is a bijection. Hence $\pi$ is an algebra 
isomorphism from $Z(R)$ onto its image $W(F,H)$ in $R(F,H)$.
\qed

\begin{lem}
If $u^{\eta}=u$, then we have $u^{\eta}v^{\eta}=(uv)^{\eta}$ for any $v\in R$.
\end{lem}
{\bf Proof:}$uv-u^{\eta}v^{\eta}=(u-u^{\eta})v+u^{\eta}(v-v^{\eta})=u^{\eta}(v-v^{\eta})\in
  RR_{\eta}(E)$. So we have $u^{\eta}v^{\eta}=(uv)^{\eta}$ for any $v\in R$.
\qed

Let $\tilde{A}$ be the subalgebra of $R$ generated by $H$. Now we
introduce a new gradation on $R$ by setting $deg(H)=1,
deg(E)=deg(F)=deg(f)+1$. This gradation is motivated by the so called
$x_{0}-$gradation suggested by Kazhdan (see \cite{K}) for the
universal eveloping algebras $U(\frak g)$ of semisimple Lie algebras
$\frak g$. Let us denote $deg(f)$ by $d$. We can define a filtration 
of $R(F,H)$ as follows:
\[
R_{(k)}(F,H)=\oplus_{i(d+1)+j\leq k} R_{i,j}(F,H)
\]
where $R_{i,j}(F,H)$ is the vector space spanned by
$F^{i}H^{j}$. Since $W(F,H)$ is a 
subalgebra of $R(F,H)$, it inherits a filtration from $R(F,H)$.
In addition, $\tilde{A}$ has the natural grading with $deg(H)=1$. 
Let us put $W(F,H)_{q}=\C-spann\{1, \Omega^{\eta},\cdots,
(\Omega^{\eta})^{q}\}$, then we have the following:
\begin{thm}
$R(F,H)$ is free as a right module over $W(F,H)$. And the
multiplication induces an isomorphism 
\[
\tilde{A}\otimes W(F,H)\longrightarrow R(F,H)
\]
as right $W(F,H)-$modules. In particular, we have the following
\[
\oplus _{p+q(d+1)=k} \tilde {A}_{p} \otimes W(F,H)_{q} \cong R(F,H)_{(k)} 
\]
\end{thm}
{\bf Proof:} First of all, $\tilde{A}\times W(F,H) \longrightarrow
R(F,H)$ is bilinear. So by the universal property of the tensor product, 
there is a linear map from $\tilde{A}\otimes W(F,H)$ into $R(F,H)$ defined by 
multiplication. It is easy to see this map is a homomorphism of right 
$W(F,H)-$modules and surjective as well. Now we show that it is
injective. Suppose $(\sum_{i=0}^{n} a_{i}H^{i})(\sum_{j=0}^{m}
b_{j}(2\eta(E)F+u(H+1))^{j})=0$ with $b_{m}\neq 0$. Then we 
have $(\sum_{i=0}^{n}a_{i}H^{i})b_{m}(2\eta(E))^{m}F^{m}+g(H,F)=0$, 
where the $F-$degree of $g$ is less than $m$ by direct computations. 
But $H^{i}F^{j}$ are part of the basis of $R(H,F)$ as a vector space, 
hence $\sum_{i=0}^{n}a_{n}H^{i}=0$. Thus, the theorem has been proved.
\qed

Let $Y_{\eta}$ be the left $R-$module defined by $Y_{\eta}=R \otimes
\C_{\eta}$ where $\mathbb{C}_{\eta}$ is the $1-$dimensional $R(E)-$module 
defined by $\eta$. It is easy to see that $Y_{\eta}=R/RR_{\eta}(E)$ is 
a Whittaker module with a cyclic vector $1_{\eta}=1\otimes 1$. And we have a 
quotient map 
\[
R \longrightarrow Y_{\eta}
\]
If $u\in R$, then $u^{\eta}$ is the unique element in $R(F,H)$ such that $u1_{\eta}=u^{\eta}1_{\eta}$.
As in \cite{K}, we define the $\eta -$reduced action of $R(E)$ on $R(F,H)$
such that $R(F,H)$ is an $R(E)$ module in the following way:
\[
x\bullet v=(xv)^{\eta}-\eta(x)v
\]
where $x\in R(E)$, $v\in R(F,H)$.
\begin{lem}
Let $u\in R(F,H)$ and $x\in R(E)$, then we have
\[
x\bullet u^{\eta}=[x,u]^{\eta}
\]
\end{lem}
{\bf Proof:} $[x,u]1_{\eta}=(xu-ux)1_{\eta}=(xu-\eta(x)u)1_{\eta}$. Hence
\[[x,u]^{\eta}=(xu)^{\eta}-\eta(x)u^{\eta}=(xu^{\eta})^{\eta}-\eta(x)u^{\eta}=x\bullet
u^{\eta}\].
\qed

\begin{lem}
Let $x\in R(E)$, $u\in R(F,H)$, and $v\in W(E,F)$, then we have 
\[
x \bullet (uv)=(x\bullet u)v.
\]
\end{lem}
{\bf Proof:} Let $v=w^{\eta}$ for some $w\in Z(R)$, then $uv=uw^{\eta}=u^{\eta}w^{\eta}=(uw)^{\eta}$.
Thus
$x\bullet(uv)=x\bullet(uw)^{\eta}=[x,uw]^{\eta}=([x,u]w)^{\eta}=[x,u]^{\eta}w^{\eta}=[x,u]v=(x\bullet u^{\eta})v=(x\bullet u)v$.

\qed

Let $V$ be an $R-$module and $R_{V}$ be the annihilator of $V$ in
$R$. Then $R_{V}$ defines a central ideal $Z_{V}\subset Z$ by 
setting $Z_{V}=R_{V}\cap Z$. Suppose that $V$ is a Whittaker module 
with a cyclic Whittaker vector $w$, we denote by $R_{w}$ the
annihilator of $w$ in $R$. It is obvious that
$RR_{\eta}(E)+RZ_{V}\subset R_{w}$. In the next
theorem, we show that the reversed inclusion holds.

First of all, we need a lemma:
\begin{lem}
Let $X=\{v\in R(F,H)\mid (x \bullet v)w=0,x\in R(E)\}$. Then
\[
X=(\tilde{A}\otimes W_{V}(F, H))+W(F, H)
\]
where $W_{V}(F,H)=(Z_{V})^{\eta}$. In fact, $R_{v}(F,H) \subset X$ and 
\[
R_{w}(F, H)=\tilde{A} \otimes W_{V}(F,H)
\]
where $R_{w}(F,H)=R_{w} \cap R(F,H)$.
\end{lem}
{\bf Proof:} Let us denote by $Y=\tilde{A} \otimes
W_{V}(F,H)+W(F,H)$. Thus we need to verify $X=Y$.
Let $v\in W(F,H)$, then $v=u^{\eta}$ for some $u\in Z(R)$. 
So $x \bullet v=x\bullet u^{\eta}=(xu^{\eta})^{\eta}-
\eta(x)u^{\eta}=(x^{\eta}u^{\eta})^{\eta}-\eta(x)u^{\eta}=((xu)^{\eta})^{\eta}-\eta(x)u^{\eta}=(xu)^{\eta}-\eta(x)u^{\eta}=x^{\eta}u^{\eta}-\eta(x)u^{\eta}=0$. So $W(F,H)\subset X$.
Let $u\in Z_{V}$ and $v\in R(F,H)$, then for any $x\in R(E)$ we have
$x\bullet(vu^{\eta})=(x\bullet v)u^{\eta}$. Since $u\in Z_{V}$, then
$u^{\eta}\in R_{w}$. Thus we have $vu^{\eta}\in X$, hence
$\tilde{A}\otimes R_{V}(F,H)\subset X$. This proves $Y\subset X$. Let
$\tilde{A_{i}}$ be the one dimensional subspace of $R(H)$ spanned by
$H^{i}$ and $\bar{W_{V}(F,H)}$ be the complement of $W_{V}(F,H)$ in
$W(F,H)$. Set $M_{i}=\tilde{A_{i}}\otimes \bar{W_{V}(F,H)}$, then we
have the following:
\[
R(F,H)=M\oplus Y
\]
where $M=\sum_{i\geq 1}M_{i}$. We show that $M\cap X= 0$. Let
$M_{[k]}=\sum_{1\leq i\leq k}M_{i}$, then $M_{[k]}$ are a filtration of
$M$. Suppose $k$ is the smallest integer such that $X\cap M_{[k]}\neq 0$ 
and $0\neq y\in X\cap M_{[k]}$. Then we have $y=\sum_{1\leq i\leq
  k}y_{i}$ where $y_{i}\in \tilde{A_{i}}\otimes \bar{W_{V}(F,
  H)}$. Now we have $0\neq E\bullet y_{i}\in \tilde{A_{i-1}}\otimes
W_{V}(F,H)$ for $i\geq 1$. Hence we have $E\bullet y \in M_{[k-1]}$. 
This is a contradiction. So we have $X\cap M=0$. Now we prove that 
$R_{w}(F, H) \subset X$. Let $u \in R_{w}(F, H)$ and $x \in R(E)$,
then we have $xuw=0$ and $uxw=\eta{x}uw=0$. Thus $[x,u]\in
R_{w}(F,H)$, hence $[x,u]^{\eta}\in R_{w}(F,H)$. Since $u\in
R_{w}(F,H)\subset R_{w}(E,F,H)$, then $x\bullet u=[x,u]^{\eta}$. 
Thus $x\bullet u \in R_{w}(F,H)$. So $u\in X$ by the definition of
$X$. Now we prove the following:
\[
W(F,H)\cap R_{w}(F,H)=W_{V}(F,H)
\] 
In fact, $W_{V}(F,H)=(Z_{V}^{\eta})$ and $W_{V}(F,H)\subset
R_{w}(F,H)$. So if $v\in W_{w}(F,H)\cap R_{w}(F,H)$, then we can
uniquely write $v=u^{\eta}$ for $u\in Z(R)$. Then $vw=0$ implies $uw=0$
and hence $u\in Z(R)\cap R_{w}(F, H)$. Since $V$ is generated
cyclically by $w$, we have proved the above statement.

Obviously, we have $R(E,F,H)Z_{V}\subset R_{w}(E,F,H)$. Thus we 
have $\tilde{A} \otimes W_{V}(F,H) \subset R_{w}(F,H)$, hence we 
have $R_{w}(F,H)=\tilde{A}\otimes W_{V}(F,H)$.
\qed

\begin{thm}
Let $V$ be a Whittaker module admitting a cyclic Whittaker vector $w$, 
then we have 
\[
R_{w}=RZ_{V}+RR_{\eta}(E).
\]
\end{thm}
{\bf Proof:} It is obvious that $RZ_{V}+RR_{\eta}(E)\subset
R_{w}(E,F,H)$. Let $u\in R_{w}(E,F,H)$, we show that 
$u\in R(E,F,H)Z_{V}+R(E,F,H)R_{\eta}(E)$. Let $v=u^{\eta}$, then 
it suffices to show that $v\in \tilde{A}\otimes W_{V}(F,H)$. 
But $v\in R_{w}(F,H)=\tilde{A} \otimes W_{V}(F,H)$.
\qed

\begin{thm}
Let $V$ be any Whittaker module for $R$, then the correspondence
\[
V\longrightarrow Z_{V}
\]
sets up a bijection between the set of all equivalence classes 
of Whittaker modules and the set of ideals of $Z(R)$.
\end{thm}
{\bf Proof:} Let $V_{i}, i=1, 2$ be two Whittaker modules. If
$Z_{V_{1}}=Z_{V_{2}}$, then clearly $V_{1}$ is equivlent to $V_{2}$ by 
the above Theorem. Now let $Z_{\ast}$ be an ideal of $Z(R)$ and let
$L=RZ_{\ast}+RR_{\eta}(E)$. Then $V=R/L$ is a Whittaker module with a 
cyclic Whittaker vector $w=\bar{1}$. Obviously, we have $R_{w}=L$. So 
that $L=R_{w}=RZ_{V}+RR_{\eta}(E)$. This implies that 
$\eta(L)=(Z_{\ast})^{\eta}=\eta (R_{w})=(Z_{V})^{\eta}$. Since $\eta$ 
is injective, thus $Z_{V}=Z_{\ast}$. Hence we have completed the proof.
\qed

\begin{thm}
Let $V$ be an $R-$module. Then $V$ is a Whittaker module if and only if 
\[
V\cong R\otimes_{Z\otimes R(E)}(Z/Z_{\ast})_{\eta}
\]
In particular, in such a case the ideal $Z_{\ast}$ is uniquely
determined to be $Z_{V}$.
\end{thm}
{\bf Proof:} If $1_{\ast}$ is the image of $1$ in $Z/Z_{\ast}$, then 
$Ann_{Z(R)\otimes R(F)}(1_{\ast})=R(E)Z_{\ast}+Z(R)R_{\eta}(E)$. The 
annihilator of $w=1\otimes 1_{\ast}$ is $R_{w}=R(E,F,H)Z_{\ast}+R(E,F,H)R_{\eta}(E)$. Then the result follows from the last theorem. 

\qed

\begin{thm}
Let $V$ be an $R-$module with a cyclic Whittaker vector $w\in V$. Then 
any $v\in V$ is a Whittaker vector if and only if $v=uw$ for some $w\in Z(R)$.
\end{thm}
{\bf Proof:} If $v=uw$ for some $u \in Z(R)$, then obviously that $v$
is a Whittaker vector. Now let $v=uw$ for some $u\in R$ be a Whittaker 
vector of $V$. Then $v=u^{\eta}w$ by the definition of Whittaker
module. So that we can assume that $u\in R(F,H)$. If $x\in R(E)$, we 
have $xuw=\eta(x)uw$ and $uxw=\eta(x)uw$. Thus we have $[x,u]w=0$ and 
hence $[x,u]^{\eta}w=0$. But we have $x\bullet u=[x.u]^{\eta}$. Thus
we have $u\in X$. We can now write $u=u_{1}+u_{2}$ where $u_{1}\in
R(F,H)$, and $u_{2}\in W(F,H)$. Then $u_{1}w=0$. Thus $u_{2}w=v$. 
But now $u_{2}=u_{3}^{\eta}$ where $u_{3}\in Z(R)$. So we have
$v=u_{3}w$, which proves the theorem.
\qed
 
Now let $V$ be a Whittaker module and $End_{R}(V)$ be the endomorphism
ring of $V$ as an $R-$module. Then we can define the following
homomorphism of algebras defined by the action of $Z(R)$ on $V$:
\[
\pi_{V}\colon Z\longrightarrow End_{R}(V)
\]
It is clear that $Z(R)/Z_{V}(R)=\pi_{V}(Z(R))\subset End_{R}(V)$. In
fact the next theorem says that this inclusion is equal as well.
\begin{thm}
Assume that $V$ is a whittaker module. Then $End_{R}(V)\cong Z/Z_{V}$. 
In particular $End_{R}(V)$ is commutative.
\end{thm}
{\bf Proof:} Let $w\in V$ be a cyclic Whittaker vector. If $\alpha \in
End_{R}(V)$, then we have $\alpha(w)=uw$ for some $u\in Z(R)$ by the
above Theorem. Thus we have $\alpha(vw)=vuw=uvw=u\alpha(w)$. Hence
$\alpha=\pi_{V}(u)$, which proves the theorem.
\qed

Now we are going to construct some Whittaker modules. Let $\xi \colon
Z(R)\longrightarrow \C$ be a central character. For any given central 
character $\xi$, $Z_{\xi}=Ker(\xi)\subset Z(R)$ is a maximal ideal of
$Z(R)$. Since $\C$ is an algebraically closed field, then 
$Z_{\xi}=(\Omega-a_{\xi})$. 
Given any central character $\xi$, let $\C_{\xi,\eta}$ be the one
dimensional $Z\otimes R(E)-$module defined by $uvy=\xi(u)\eta(v)y$ for
any $u\in Z$ and $v\in R(E)$. Let
\[
Y_{\xi,\eta}=R(E,F,H)\otimes_{Z\otimes R(E)} \C_{\xi,\eta}
\] 
It is easy to see that $Y_{\xi, \eta}$ is a Whittaker module of type
$\eta$ and admits an infinitesimal central character $\xi$. Since $R$
is almost commutative, so by Quillen's lemma, we know every
irreducible representation has an infinitesimal central character. As
studied in \cite{S}, we know $R$ has a similar Verma module theory. 
In fact, Verma modules also fall into the category of Whittaker
modules if we allow the trivial central character on $R(E)$. Namely, 
we have 
\[
M_{\lambda}=R\otimes_{R(E,H)} \C_{\lambda}
\]
where $R(H)$ acts on $\C_{\lambda}$ through $\lambda$ and $R(E)$ acts
trivially on $\C_{\lambda}$. $M_{\lambda}$ admits an infinitesimal
character with $\xi=\xi(\lambda)$. It is well-known that Verma modules
may not be necessarily irreducible, even though they have
infinitesimal central characters. While Whittaker modules 
are on the othe extreme as shown in the next theorem:
\begin{thm}
Let $V$ be a Whittaker module for $R$. Then the following conditions
are equivalent.
\begin{enumerate}
\item $V$ is irreducible.\\
\item $V$ admits a central character.\\
\item $Z_{V}$ is a maximal ideal.\\
\item The space of Whittaker vectors of $V$ is one-dimensional.\\
\item All nonzero Whittaker vectors of $V$ are cyclic.\\
\item The centrializer $End_{R}(V)$ is reduced to the constants $\C$.\\
\item $V$ is isomorphic to $Y_{\xi,\eta}$ for some central character $\xi$.
\end{enumerate}
\end{thm}
{\bf Proof:} It is easy to see that $(2)-(7)$ are equivalent to each other
by using the above theorems we have just proved. We also know $(1)$
implies $(2)$. To complete the proof, it suffices to show that $(4)$
implies irreducibility. To this true, we show that any submodule $V^{\prime}$
of $V$ contains a nonzero Whittaker vector, which closes the proof.
Let $v\in V$, we recall that the reduced $\eta-$action is defined as follows: 
\[
x\bullet v=xv-\eta(x)v
\]
for any $x\in R(E,F,H)$. If $u\in R$ and $x\in R(E)$, then we have
$x\bullet(uv)=xuv-\eta(x)uv=[x,u]v+uxv-\eta(x)uv$. Since
$uxv=\eta(x)uv$, thus we have the following:
\[
x\bullet(uv)=[x,u]v
\]
Now let 
\[
R(E,F,H)\longrightarrow V
\] 
be the morphism from $R(E,F,H)$ into $V$ by 
mapping $u\in R$ to $u\bullet w$.
Then this map is a homomorphism of the $R(E)-$module $R$ under the adjoint
action of $R(E)$ into the $R(E)-$module $V$ under $\eta-$reduced action.
Note the adjoint actin of $R(E)$ on $R$ is locally finite. Let $0\neq
v\in V^{\prime}$ and write $v=uw$ for $u\in R$. Let $R_{0}$ be the
$R(E)-$submodule of $R$ generated by $u$. Then the submodule $R_{0}\subset
R$ is finite dimensional. Thus the image $V_{0}$ of $R_{0}$ inside $V$ 
is finite dimensional. And $R(E)$ is the enveloping algebra of the one
dimensional Lie algebra generated by $E$, which acts nilpotently on
$V_{0}$ via the reduced action. Since we have $v\in V_{0}$, then
$V_{0}\subset V^{\prime}$. So by Engel's Theorem, we have $x\bullet
v_{0}=0$ for some $0\neq v_{0} \in V_{0}$ for all $x\in R(E)$. 
So $v_{0}$ is a Whittaker vector.
\qed

It is easy to prove the following theorems, we refer the reader to
\cite{K} for details about their proofs:
\begin{thm}
Let $V$ be an $R-$module which admits an infinitesimal
characater. Assume that $w\in V$ is a Whittaker vector. 
Then the submodule $Rw\subset V$ is irreducible.
\end{thm}

\qed

\begin{thm}
Let $V_{1},V_{2}$ be any two irreducble $R-$modules with the same
infinitesimal character. If $V_{1}$ and $V_{2}$ contain Whittaker 
vectors, then these vectors are uniquge up to scalars. Furthermore, 
$V_{1}$ and $V_{2}$ are isomorphic to each other as $R-$modules.
\end{thm}

\qed

In fact, we have the following description about the basis of any 
irreducible Whittaker module $(V,w)$, where $w\in V$ is a cyclic 
Whittaker vector. 
\begin{thm}
Let $(V,w)$ be an irreducible Whittaker module with a Whittaker vector
$w$, then $V$ has a $\C-$vector space basis consisting of elements
$\{H^{i}w \mid i\geq 0\}$.
\end{thm}
{\bf Proof:} Since $w$ is a cyclic Whittaker vector of the Whittaker 
module $(V,w)$, then we have $V=Rw$. Since 
$R=R(F,H)\otimes_{R(E)} \C_{\eta}$, then we have $V=R(F,H)w$. 
Since $(V,w)$ is irreducible, then $(V,w)$ has a central infinitesimal 
character. Thus we have $\Omega w=\lambda(\Omega)w$. Now 
$\Omega w=(2\eta(E)F+u(H+1))w$. Hence the action of $F$ on 
$(V,w)$ is uniquely determined by the action of $H$ on $(V,w)$. Thus the 
theorem follows. 

\qed

\section{The submodule structure of  a Whittaker module $(V,w)$}
In this section, we spell out the details about the structure of 
submodules of a Whittaker module $(V,w)$. We have a clean description 
of all submodules through the geometry of the affine line
${\mathbb{A}}^{1}$. Throughout this section, we will fix a Whittaker 
module $V$ of type $\eta$ and a cyclic Whittaker vector $w$ of $V$. 
Our arguement is more or less the same as the one in \cite{O}.
\begin{lem}
Let $Z(R)=\C[\Omega]$ be the center of $R$, then any maximal ideal 
of $Z(R)$ is of the form $(\Omega-a)$ for some $a\in \C$. 
\end{lem}
{\bf Proof: } This fact follows from the assumption that $\C$ is
algebraically closed field and Hilbert Nullstenllenzuts Theorem.
\qed

Let $Z_{V}$ be the annihilator of $V$ in $Z(R)$, then
$Z_{V}=(f(\Omega))$ for some polynomial $f(\Omega)\in Z(R)$. 
Suppose that $f=\prod_{i=1,2, \cdots, m} f_{i}^{n_{i}}$ for 
some irreducible polynomials $f_{i}$. Then we have the following:
\begin{prop}
$V_{i}=R\prod_{j\neq i} f_{j}^{n_{j}}w$ are indecomposable submodules 
of $V$. In particular, we have
\[
V=V_{1}\oplus \cdots \oplus V_{m}
\] 
as a direct sum of submodules.
\end{prop}
{\bf Proof:} It is easy to verify that $V_{i}$ are submodules. 
Now we show each $V_{i}$ is indecomposible. Suppose not, we can 
assume without loss of generality that $V_{1}=W_{1}\oplus W_{2}$. 
Note that $ Z_{V}=Z_{W_{1}}\cap Z_{W_{2}}$. Since $Z(R)$ is a
principal ideal domain, hence $Z_{W_{i}}=(g_{i}(\Omega)$. Thus 
we have $g_{i}\mid f_{1}^{n_{1}}$. This implies that the decomposition 
is not a direct sum. Therefore, $V_{i}$ are all indecomposable. The 
decomposition follows from the Chinese Reminder Theorem.
\qed

\begin{prop}
Let $(V,w)$ be a Whittaker module and $Z_{V}=<f^{n}>$ where $f$ is an
irreducible polynomial in $\C[\Omega]$. Let $V_{i}=Rf^{i}w, i=0,\cdots
n$ and $S_{i}=V_{i}/V_{i+1}, i=0,\cdots, n-1$. Then $S_{i},i=0,\cdots,
n-1$ are irreducible Whittaker modules of the same type $\eta$ and
form a composition series of $V$. In particular, $V$ is of finite length. 
\end{prop}
{\bf Proof:} The proof follows from the fact that $Z_{S_{i}}=<f^{i}>$ 
for all $i$.
\qed

\begin{rem}
$V$ being of finite length is an analogue of the classical situation. 
In deed, in the classical case, Whittaker modules of $U(\frak g)$ are
mapped to holonomic $D-$modules on the flag variety of $\frak g$ by
the Beilinson-Bernstein localization
(\cite{BB}), and therefore are of finite length (\cite{MS1} and \cite{MS2}).
\end{rem}

With the same assumption, we have the following
\begin{cor}
$V$ has a unique maximal submodule $V_{1}$.
\end{cor}
{\bf Proof:} This is obvious, since the only maximal ideal of $Z_{V}$
is $<f>$.
\qed

Based on the above propositions, the irreducibility and
indecomposibility are reduced to the investigation of $Z_{V}$. 
One has that $V$ is irreducible if and only if $Z_{V}$ is a maximal 
ideal. And $V$ is indecomposible if and only if $Z_{V}$ is primary. 
The following proposition is just a refinement of the submodule 
struture of $(V,w)$.
\begin{prop}
Suppose $(V,w)$ is an indecomposible Whittaker module with
$Z_{V}=<f^{n}>$, then any submodule $V^{\prime} \subset V$ is of the 
form:
\[
V^{\prime}=Rf^{i}w
\] 
for some $i\in \{0,\cdots,n\}$.
\end{prop}
\qed

Now we are going to investigate the submodule structure of any
Whittaker module $(V,w)$ with a nontrivial central annihilator
$Z_{V}$. First of all, we recall some notations from \cite{K}. 
Let $V^{\prime}\subset V$ be any submodule of $V$, we define an 
ideal of $Z$ as follows:
\[
Z(V^{\prime})=\{x\in Z\mid x V^{\prime}\subset V^{\prime}\}
\]
We may call $Z(V^{\prime})$ the normalizer of $V^{\prime}$ in $Z$. 
Conversely, for any ideal $J\subset Z$ containing $Z_{V}$, $JV\subset
V$ is a submodule of $V$.
\begin{thm}
Let $(V,w)$ be a Whittaker module with $Z_{V}\neq 0$. Then there is a 
one-to-one correspondence between the set of all submodules of $V$ and 
the set of all ideals of $Z$ containing $Z_{V}$ given by the maps
$V^{\prime}\longrightarrow Z(V^{\prime})$ and $J\longrightarrow JV$. 
These maps are inverse to each other.
\end{thm} 
{\bf Proof:} The proof is straightforward.
\qed

Now we have a description of the basis of any Whittaker module
$(V,w)$ as follows.
\begin{prop}
Let $(V,w)$ be a Whittaker module and suppose that $Z_{V}=<f(\Omega)>$ 
where $f\neq 0$ is a monic polynomial of degree $n$. Then ${\bf
  B}=\{F^{i}H^{j}w \mid 0\neq i\leq n-1, j\in \Z_{\geq 0}\}$ 
is a basis of $(V,w)$. If $f=0$, then ${\bf B}=\{F^{i}H^{j}w\mid i,
j\in \Z_{\geq 0} \}$ is a basis of $(V,w)$.
\end{prop}
{\bf Proof:} The proof is easy.

\qed

\end{document}